\documentclass{article}
\usepackage{hyperref}
\usepackage{amsmath}
\usepackage{amscd}
\usepackage{amsthm}
\usepackage{amssymb} \usepackage{latexsym}
\usepackage{eufrak}
\usepackage{euscript}
\usepackage{graphicx} 
\usepackage{array}
\usepackage{enumerate}
\usepackage{xy}
\usepackage{makeidx}
\usepackage{pictexwd,dcpic}

\usepackage{color}

\newcommand{\bel}[1]{\begin{equation}\label{#1}}

\newcommand{\be}{\begin{equation}}

\newcommand{\ba}{\begin{eqnarray}}
\newcommand{\ea}{\end{eqnarray}}
\newcommand{\rf}[1]{(\ref{#1})}
\newcommand{\bi}{\bibitem}
\newcommand{\qe}{\end{equation}}

\title{Object oriented models vs. data analysis -- is this the right
  alternative?}
\author{J\"urgen Jost}
\begin{document}
\maketitle

\section{Introduction: The basic issue}
Traditionally, there has been the distinction between pure and applied mathematics. Pure mathematics -- so the story goes -- discovers, creates and investigates abstract structures (see \cite{J8}) for their own sake, while applied mathematics applies existing mathematical tools and develops new ones  for specific problems arising in other sciences. The interaction between pure and applied mathematics takes place in both directions. Applied mathematics utilizes concepts and methods developed in pure mathematics, and problems from diverse applications in turn stimulate the development of new mathematical theories. And traditionally, those problems arose within a clear conceptual framework of a particular science, most notably physics. In this essay, I want to argue that this distinction between pure and applied mathematics is no longer useful -- if it ever was --, and that the challenge of  large and typically rather diverse data sets, typically arising from new technologies instead of   theoretically understood and experimentally testable concepts, not only calls for new mathematical tools, but also necessitates a rethinking of the role of mathematics itself. 

In fact, already in the past, what was called applied mathematics often  developed domain independent methods. And this domain independence typically led to a gain in generality and mathematical depth. Statistics, for instance, has become so powerful precisely because it developed methods and concepts that apply to essentially any field, from particle physics to the  social sciences. And the error estimates and convergence rates  provided by numerical analysis are valid for any application of a particular numerical method in the engineering sciences or elsewhere. Nevertheless, such applications usually took place within a particular field with a well developed theoretical framework that provided interpretations for the statistical or computational results obtained. And in other cases, like modeling with differential equations, the specific properties of a concrete scientific theory yielded key ingredients for the mathematical structures. The development of the underlying scientific theory and the mathematical model often went hand in hand, and they deeply depended upon each other. Mathematical modeling was an essential ingredient of the scientific strategy, and  in that sense, mathematics was more than a tool. In contrast, nowadays mathematicians are  often confronted with data sets of obscure quality, perhaps even of dubious origin, and without any firm theoretical foundation. Even the distinction between meaningful data and meaningless noise may not be clear at all.

Therefore, as data collection these days is typically ahead of theoretical understanding, mathematics  should radically face the lack of theory and look at what there is, the data, and see what it can do with them. And what I want to call for are not ad hoc methods for every concrete data set, but rather an abstract analysis of the deeper structural challenges. Of course, modern mathematics is developed and sophisticated enough to provide appropriate and helpful tools for basically any data set, but this by itself is too narrow a scientific perspective.  For me as a mathematician, mathematics is more than data analysis. We need a conceptual rethinking.

Let us take a look at the situation from the perspective of the sciences. It seems to me that in the presence of data, there are two different, but interwoven issues:
\begin{enumerate}
\item The epistemic question, or the role of theory: Can we have data without a theory? And if
  not, does the theory have to be specific for the domain from which
  the data are collected?
\item The ontological question, or the role of models: Do we need, or have to postulate, specific
  objects underlying the data? 
What is it that the data tell us something about, and how can we, or
should we model that? 
\end{enumerate}
This distinction may sound a little like issues debated between
scientific realists, positivists, and empiricists. But even if we adopt the
latter stance and use, for instance, van Fraassen's criterion that a
``theory is empirically adequate if it has some model such that 
all appearances are isomorphic to empirical substructures of that
model'' where ``appearances'' are ``structures which can be described
in experimental and measurement reports'' (\cite{vF}, p.64), the
problem of data without such an empirically adequate theory  remains. 

Of course, there are some easy answers that seem rather obvious:
\begin{itemize}
\item Data that cannot be interpreted within some theoretical
  framework are meaningless. Much of current data collection is solely
  driven by particular technologies rather than by scientific
  questions. As Sir Peter Medawar put it, ``No new principle will
  declare itself from below a heap of facts''.\footnote{I learned this quote from Peter Schuster.} And for ``facts'', we might want to substitute ``data''.
\item If we simply analyze data with intrinsic formal tools whose
  precise functioning we may not even understand, then science becomes
  agnostic \cite{Nap}. We need to know first about what the data are telling us
  something. That is, we need an underlying structure from which the
  phenomena revealed by the data are derived. The meaning of data
  depends on the context from which they are collected. 
\end{itemize}

These concerns are clearly valid. Nevertheless, I think that the
issues are somewhat more subtle.
\begin{enumerate}
\item To what extent does theory have to be domain specific? Or more
  precisely, what is the relationship between general theoretical
  issues -- to be elaborated below -- and domain specific ones? Or
  more positively, will the current situation of ``big data'' lead to
  new types of theories that start from a data intrinsic rather than a
  domain specific perspective, and could that possibly lead to
  theoretical insights at a higher level of abstraction?
\item What is the right balance between a purely phenomenological
  approach and one that starts with a model involving underlying
  objects? Such objects could be either  the carriers of the properties revealed by the data or at least 
  offer  formal structures to which quantities measured on the data
  are isomorphic. 
\end{enumerate}

A traditional
perspective would consider mathematics only as a tool when -- possibly
very large -- data sets are to be analyzed. One of the theses of this
essay is that this process itself, data analysis, can and already has become an object
of mathematical research. Thus, mathematics not only serves as a
powerful tool, but by reflecting this role, gains a new perspective
and lifts itself to a higher level of abstraction. Then, the domain of
such mathematical inquiry no longer is a specific field, like physics,
but mathematical technique itself. Mathematics then is no longer, if it ever was, a mere formal tool for science, but becomes the science of formal tools. 

Another thesis of this essay is
more concerned with the traditional role of mathematics as the formal
analysis of models arising from specific domains. I shall
argue here against the principle of reductionism. The thesis will be
that in the empirical domains where regularities have been identified
that could be cast into mathematical structures, the relations between 
the different levels at which such structures emerge are very
intricate. A simple formal structure at one level typically can
neither be derived from an even simpler structure at a lower level,
nor can be used to identify the relevant structure at a higher
level. The structures at both lower and higher levels can be
substantially more complex. Again, it is a general task for
mathematical research to identify and more deeply understand the
principles of such transitions between different scales and
levels. That is, what is at issue are not so much the phenomena at the different levels, but rather the transition laws between levels. Much is known here already, both from the side of mathematics
and that of physics, but a more abstract and fundamental mathematical
theory is still missing. 

There might even be unexpected connections between those two aspects,
a mathematics of and not only for data analysis on one side, and the principles of
transitions between different levels on the other side. At least,
there is the question of how such transitions can be inferred from the
data. I regard this as a fundamental question for any abstract theory
of data. 

Before addressing these issues in detail, I should point out that there is another, perhaps even more fundamental and important role for mathematics in the sciences, although that role will not be discussed in this essay. This consists in providing a framework for conceptual thinking and formal reasoning. This is valid across all disciplines, and actually badly needed in many of them. I shall address this aspect elsewhere.

\section{The role model of  physics}\label{physics}
In some way or another, there is a particular role model behind many
discussions, that of classical physics. According to that perspective, every science should strive
towards such a model. We should have  objects whose properties and interactions can
be formally described by mathematical relations, in particular,
differential equations. And we should be able to conduct experiments
whose outcomes can be understood within an encompassing theoretical
framework. Deviations from that role model are regarded as deficits.

Newton's theory combined the terrestial dynamics as developed by
Galilei and others and Kepler's laws of celestial motions in a unified
theory of gravity. The basic objects were pointlike masses. Such a
concept would have appeared meaningless to Descartes who considered
extension as the basic quality of matter. In Newtonian dynamics,
however, the essential feature is the mass of an object, and this
allowed him to work with extensionless points carrying a mass as
idealized objects. Newton's theory allows for an exact solution of the
two-body problem.\footnote{The one-body problem, however, was eventually
understood as being ill posed. As Leibniz first saw, any physical theory
has to be concerned with relations between objects, and an irreducible
point mass can only entertain relations with other such objects, but
not with itself, because there is no fixed external frame of reference
independent of objects. The latter aspect is fundamental in Einstein's theory
of general relativity. (See for instance the exposition by the author
in \cite{Rie}.)} Besides points, there are other
fundamental constituents of physical theories. These are the fields,
like the electromagnetic one. In physical theories, fields are
considered to be no less real than material objects. In particular,
they also carry forces which make physical interactions possible. 
Also, while the two-body problem admits an exact solution in Newton's
theory, this is no longer true for the three-body problem. Although
that problem is perfectly well posed within Newton's theory of
gravity, it can in general no longer be solved in closed form, unless
particular symmetries pertain. This was first realized  by
Poincar\'e. Nevertheless, one can develop approximation schemes or
compute numerical solutions of the differential equations to any
desired degree of accuracy with sufficiently large computer power. At
a more theoretical level, the stability against perturbations becomes
a subtle issue, as developed in KAM theory (after Kolmogorov, Arnol'd,
and Moser), see \cite{SM,TZ}. Also, there are mathematical examples \cite{Xia,SX} where internal
fluctuations in such a system can grow without bounds until one of the
bodies is ejected to infinity. 

The paradigmatic objects of this theory are the celestial bodies in
the solar system, and the theory is concerned with their orbits. Even when the theory is restricted to the basic
configuration of the sun, the earth, and the moon, there are more than
two such objects. And these objects are by no means pointlike, but
rather extended and possessing a non-homogeneous internal
structure. If one models them in more detail, the resulting theory
becomes much more complicated, and at best one may hope for numerical
approximations, but these will inevitably be of limited scope and
validity. 

But this then raises the question why such a simple theory as Newton's
is applicable at all to such a complicated configuration, and in
particular, why numerical solutions of Newton's equations lead to such
accurate descriptions and predictions of the movements of the
planets. (For the sake of the argument, we ignore here the corrections
necessitated by Einstein's theory.) When we change the scale, either
going down to the internal structure of the celestial bodies, or up to
configurations of many gravitating objects, no such simple theory
applies. 

Is Newton's theory then a completely singular instance, and would we
therefore be ill-advised to consider it as a role model of a
scientific theory? Well, it is not completely singular, and in order
to understand the issues involved better, let us consider another
example. 

The fundamental equation of quantum mechanics is Schr\"odinger's
equation
\begin{equation}
\label{fi1}
\sqrt{-1} \hbar \frac{\partial \phi(x,t)}{\partial t}=-\frac{\hbar^2}{2m}\Delta \phi(x,t)+V(x)\phi(x,t)
\end{equation}
for the quantum mechanical state $\phi(x,t)$ at position $x$ and time
$t$, where $\Delta$ is the Laplace operator, a second order partial
differential operator, $V(x)$ is the potential at $x$, and $\hbar,
m$ are physical constants.
The Schr\"odinger equation no longer describes the deterministic behavior of
pointlike objects, but rather the evolution of probabilities. The
dynamics of these probabilities are still deterministic. And in fact,
for the hydrogen atom, the Schr\"odinger equation is exactly
solvable. For more complicated atoms, we encounter the difficult
numerical problems of quantum chemistry. The reason why it applies so
well to the hydrogen atom is that its nucleus consists of a single
proton which under normal conditions essentially behaves like an 
elementary, indivisible (``atomic'' in the literal sense) particle. In
fact, however, the proton is not elementary, but is composed of more
elementary particles, the quarks. That is, the fact that the hydrogen
atom can be described by the Schr\"odinger equation does not have a
theoretical, but an empirical reason, the so-called quark
confinement. Of course, one may then try to derive this from the
standard model of particle physics, but that is a different
issue. (For the theoretical background, see e.g. \cite{Wein}.) What remains is the empirical validity of  Schr\"odinger's
equation. The standard model, in contrast, is a model whose deeper
justification still is a matter of intense debate. There are
approaches like string theory (see e.g. \cite{GSW,Pol,J5,J4,Zwie})  which, however, are beyond the range of
experimental testing for the time being. Why Schr\"odinger's
equation applies so well (at least to the hydrogen atom) then remains
mysterious at a deeper level. And it even applies to other atoms,
although it is no longer exactly solvable already for the Helium atom.
This is analogous to Newton's equation and the three-body
problem. 

The equations of Newton and Schr\"odinger may possess certain universal properties, and some aspects might be  considered with the modern tools of renormalization groups, see for instance \cite{Bat}, but the question why such equations appear at particular scales and not at others remains. Thus, my question is not why a mathematical description of nature is possible, or why mathematics is so effective as a tool for describing physical reality, but rather why it works so well at particular scales or under certain circumstances, but not in other situations. Historically, when the mechanical natural philosophy of the 17th century could capture so well the movements of celestial bodies, explain the acceleration of falling bodies or the parabolic trajectories of canon balls or identify the conservation laws of inelastic collisions, it was expected that this success could easily be extended to other domains and achieve, for instance, a mechanical explanation of animals, and, as some, like La Mettrie, believed, even of humans. Of course, this scientific program turned out to be a miserable failure, and biology as a modern science got off the ground only in the 19th century, on the basis of completely different principles than those of the mechanical natural philosophy of the 17th century. \\

I shall return in Section \ref{valcomp} to the difference between the theoretical validity of
an equation and its exact solvability, as this is important for our
topic. Before that, however, I'll proceed to larger physical
scales. In the area of the physics of nano- and microstructures, on
the one hand, there are atomic or molecular models, and on the other hand,
there are continuum models. Either type of model can yield accurate
descriptions within its range of validity and make precise
predictions. One might then argue that the discrete (atomic or
molecular) models are the basic ones that capture the underlying
objects whereas the continuum models are purely
phenomenological. This, however, will miss the fundamental
mathematical question, the relationship between these two types of
models and the transition between them when changing the scale. Again,
this is an aspect to which I need to return. For instance, it seems that
in many respects, the Navier-Stokes equations correctly -- whatever
that means -- or at least adequately describe the behavior of
fluids. These are continuum equations. The numerical solution, that is, an approximation by a numerical scheme, is one of the most important topics of scientific computation, and by now, there exist numerical schemes that can solve them to a very high degree of accuracy. Nevertheless, within the mathematical theory of partial differential equations, the general existence of a solution to date has not been established, and in fact, this is considered to be a challenging and difficult problem. The issue of turbulence which
is still not understood reveals that there are some deep
problems lurking here. Interestingly, while turbulent flow seems to be
rather chaotic, there might exist some statistical regularities. In any case, the Navier-Stokes equation serve as the most important model in fluid dynamics.

When we move further up with regard to the physical scale, we come to
the example already discussed, the dynamics of the orbits of the
planets in the solar system which (except for relativistic
corrections) are so well described by Newtonian mechanics. Planets,
however, are no elementary bodies, but composite objects. So, let us
repeat the question why such complicated objects as planets, just
think of the geology of the earth, follow such elementary laws. If we
go to a smaller scale,  everything breaks down, or at least becomes
substantially more complicated. For instance, which mathematical or
physical theory can explain the rings of Saturn? 

And if we go to configurations of large numbers of gravitating bodies,
Newton's theory, although applicable in principle (modulo relativistic
corrections, again), becomes useless. We should need to change to a
statistical description \`a la Boltzmann \cite{Bol}. Such a statistical
description would still be mathematical; it would derive its power
from ignoring the details of the lower level structure. 

Physics then loses some of its special role. It seems to be a general
phenomenon that in many domains   particular scales exist at which
a rather simple mathematical model can capture the essential aspects,
sometimes 
even   accurately in quantitative terms, while that need no longer hold for 
smaller or larger scales. It might seem that this occurs less often in
domains other than physics, but there do exist positive examples. In
chemistry, we have the reaction kinetics. For instance, reaction
kinetics of Michaelis-Menten type are important in biochemistry. They depend on important simplifications. In  particular, they assume that the chemical substances participating in the process in question are uniformly distributed across the cell. The internal spatial structure of the cell is completely ignored. 
As another example, the
differential equations of Hodgkin-Huxley type model \cite{HH} the generation of
the spike of a neuron in a quantitatively exact manner, see
e.g. \cite{Mur,Koch}. The model operates at the cellular level and describes the dynamics there in a deterministic manner. If one goes
down to the lower, molecular, scale, things get much more complicated,
and the deterministic quantities in the Hodgkin-Huxley equations
become stochastic, that is, probabilities for the opening and closing
of certain ion channels. Similarly, if one moves up and considers
interacting systems of neurons, things become quite complicated, and
completely different types of models are called for. Thus, we have a  relatively simple, formally closed, and quantitatively accurate model at the cellular level, but no such model at the smaller molecular or the larger tissue scale.  
I shall present the mathematical structure of biochemical reaction kinetics and the biophysics of the Hodgkin-Huxley model in detail below -- and readers not interested in those mathematical or physical details can skip those --, but let me first emphasize again the general point. Not only in physical, but also in biological systems, there may exist one or more particular scales at which a gross simplification can lead to simple, but nevertheless quantitatively accurate models, and these scales then correspond to particularly useful levels of description. Typically, however, this is no longer so when we move to either smaller or larger scales. This leads to the question why and how such particular levels emerge, at which such simple quantitative models work. And also, why are such levels apparently so rare?

Let
us now describe the mechanism of biological reaction kinetics in some detail to see what is
involved. General references for biochemical kinetics are \cite{Mur},
\cite{Klipp}, and we follow here the presentation in \cite{J6}. The
basis  is the law of mass action  which states that the reaction rate
of a chemical reaction is proportional to the concentrations of the
reactants raised to the number in which they enter the reaction. That
expression is proportional to the  probability that the reactants
encounter and react with each other. Let us consider  the simple reaction
\begin{equation}
\label{kin1}
S_1+S_2 \leftrightharpoons 2P
\end{equation}
that converts
$S_1+S_2$ into $2P$ with forward rate  $k_+$ and backward rate
$k_-$. That is, when a molecule of substance $S_1$ encounters one of
$S_2$, they react and form two copies of $P$ with rate
$k_+$. Conversely, two $P$-molecules together can decay into a copy of
$S_1$ and a copy of $S_2$ at the rate $k_-$. Thus, the chemical
reaction can take place in either direction, but at possibly different
rates. If we denote the respective concentrations by $s_1,s_2,p$ and
consider them as functions of time, then
\rf{kin1} leads to the differential equations
\ba
\nonumber
\frac{ds_1}{dt}&= -k_+ s_1 s_2 +k_- p^2\\
\nonumber
\frac{ds_2}{dt}&= -k_+ s_1 s_2 +k_- p^2\\
\label{kin2}
\frac{dp}{dt} &=2(k_+ s_1 s_2 -k_- p^2)
\ea
whose solutions $s_1(t),s_2(t), p(t)$ then give the respective
concentrations at time $t$ as functions of their initial values. 
In enzymatic reactions, there also is the complex $ES$ of the enzyme
$E$ and the substrate $S$. The Michaelis-Menten theory makes the
simplifying 
assumption of a quasi-steady state for the complex $ES$, that is, its
concentration is not changing in time, and this concentration can then
be simply computed. This assumption, that the concentration of an
intermediate product remains constant, reduces a multidimensional process to a single
equation. The resulting simple systems of ordinary differential
equations capture the concentration level of substances involved in
biochemical reactions in the cell well enough for many purposes. In
contrast to the Schr\"odinger equation which considers the state
$\phi(x,t)$ as a function of the position $x$ and time $t$, the
chemical concentrations $s_1,s_2,p$ in \rf{kin2} are considered as
functions of $t$ only. That is, the model assumes that they are
homogeneously distributed in the cell. This means that  the concentrations are
assumed to be the same across the cell at each fixed time $t$. This
is, of course, a gross simplification, and the model can be refined by
allowing for varying concentrations and diffusion effects. The
point I want to make here, however,  is that even under this
 (and the further Michaelis-Menten type) simplification,  the
 solutions of  systems like \rf{kin2} can effectively
describe the values of concentrations of chemical reactants in cells.

 We next turn to the other example discussed above, the Hodgkin-Huxley model \cite{HH} for the generation of spikes in neurons,\footnote{The original model was developed only for a particulat type of neuron, the giant squid axon, but similarly models have subsequently been developed for other classes of neurons as well.} and  give a brief description of the model and
its dynamics, see e.g. \cite{J6,J7} for more details.  Formally, the
model consists of four coupled (partial) differential equations, but
they are more difficult than those for biochemical reactions, and so,
we refrain from writing them down here. Instead, we only discuss the
qualitative features of the model. The basic dynamic variable is the
membrane potential $V$ of the neuron. The potential is caused by  the different densities of charged ions in-
and outside the cell. The boundary of the cell, the cell membrane, is
impermeable to most charged ions, except for channels that are
selectively permeable for certain specific ions. That permeability in
turn depends on the membrane potential as well as on the concentration of
certain intracellular and extracellular substances. 
The differential equations of Hodgkin and Huxley  then link the temporal changes of these
quantities,  that is, the membrane potential and three gating variables
that control the permeability of the cell membrane for specific electrically charged ions. Concerning the evolution of the membrane potential,
the idea is that the time derivative of this potential is proportional
to the derivative of the electrical charge, hence to the current
flowing. That current in turn is the sum of an internal membrane current
and an external current. The latter  represents the external input to the cell. The internal 
membrane current is a sum of specific terms, each proportional to the
difference between the potential and some specific rest term. These
proportionality factors then depend on the corresponding gating
variables. Thus, the dynamical interplay between the membrane current and the gating variables is the key point. When the system is near its resting value and some positive current is injected that lifts the membrane potential $V$ above some threshold, then a positive feedback between $V$ and the fast one, $m$,  among the gating variables sets in. That is, the potential $V$ rises, and positively charged sodium ions flow in, rising the potential further. Soon, the neuron emits a spike. But then, the two slower gating variables, $h$ and $n$, take over, causing an inactivation of the inflow of the sodium ions and reversing the potential by an outflow of potassium ions. The potential then drops even below its resting value, but after some time, during which no further spike is possible, recovers to that latter value, and the neuron is receptive again to repeat the process in response to some new external current. The interaction of two time scales is important for these dynamical properties. The positive feedback between the fast variables $V$ and $m$ triggers the spike, whereas the slow variable $h$ ultimately stops the inflow of positive ions, and $n$ causes other positive ions to flow outwards, to reverse the effect and bring the potential back to (or, more precisely, below) its resting value. As already mentioned, this particular model was conceived for the giant squid axon. For this neuron, quantitative measurements were easier than for other neurons, and therefore, Hodgkin and Huxley could carefully fit the parameters of their model. Following the pioneering work of Hodgkin and Huxley, then also models for other classes of neurons were developed. While the details are different, the essential principles are the same. Thus, we have quantitative accurate models of the biophysics of neurons that operate at the cellular level, even though the details at the smaller, molecular level are much more intricate. No such simple models exist at that scale.

\section{Validity vs. computability}\label{valcomp}
We return to Poincar\'e's insight into the impossibility of an exact
solution (a solution in closed form, that is, in the form of an
explicit formula for the positions of the bodies involved at every
instance of time) of the three-body problem. Even if a mathematical
theory of nature were exactly valid in principle -- let us ignore the
quarks and the internal structure of the planets for the sake of the
argument here and assume that for instance, Newton's or
Schr\"odinger's 
equations were exactly valid  --, it could not offer an exactly
solvable description of the dynamical behavior of its objects. It is a
mathematical question when and under which particular conditions an 
exact solution by an explicit formula is possible. This has nothing to
do with the range of validity of the theory in question. Almost all
differential equations cannot be solved by an explicit formula (they
do not constitute completely integrable systems in the language of
classical mathematics). The search for closed solutions, that is, to
demonstrate the complete integrability of a dynamical systems, was one
of the highlights of 19th century mathematics. In the 20th century, in
the wake of Hilbert, mathematics turned to a different
approach. Instead of attempting to construct an explicit solution of a
specific differential equation, the problem was converted into showing
the existence of solutions for large classes of differential
equations. This was the guiding theme in particular for partial
differential equations, see for instance \cite{J1}. A key point was to
separate the abstract question of the existence from the explicit
representation of a solution. 
Once the existence of a solution had been demonstrated from abstract
principles, it then became the task of mathematics to understand the
properties of such abstract solutions, and in particular to find out
under which conditions singularities can be avoided, and then to derive
approximation schemes for such solutions and to convert them into
algorithms for their numerical construction with a precise error
control. This is basic for all the modern applications of partial
differential equations in the engineering sciences and elsewhere. This does not mean, however, that mathematics becomes purely
instrumental. 

A similar issue arises for the weather forecast. Again, it was
Poincar\'e who first realized the specific aspects of chaotic
dynamics. Some meteorological models like the Lorenz equations indeed
exhibit chaotic behavior. Therefore, we cannot accurately predict next
week's weather, because of the amplifications of tiny fluctuations by
chaotic dynamics. It is claimed, however, that the global climate 50
years from now can be predicted within certain bounds, for any
scenario of greenhouse gas emissions. The reason why this is feasible is that on a longer time scale, the daily weather fluctuations average out, and only the long term trend remains, and that is precisely what the climate models try to capture. In any case, the pecularities of chaotic dynamics as exhibited by  weather models  are not an artefact caused by any deficiency of a model, but are -- according to deep mathematical insights -- rather grounded in specific structural features of the pertinent equations. The
mathematics of chaos goes much beyond the rather simple exponential
divergence of individual trajectories, and gains positive insight via
such concepts as invariant measures (for details, see e.g. \cite{J2}). It can even elucidate the
universal nature of many aspects of chaotic dynamics, thereby going
much beyond any specific object domain. 

In this section, I have described two different trends that distinguish 20th century from 19th century mathematics. On the one hand, the ideal of writing down an explicit solution has been abandoned for the less ambitious aim of approximating a   solution to any specified degree of accuracy, or in practice at least as accurately as the computing facilities permit. When one has a good theoretical control of the numerical scheme, then higher accuracy simply requires more computer power. In that regard, mathematics has become more instrumental, concentrating on the computational tools rather than on the explicit formulae. On the other hand, general structural insights into chaotic dynamics teach us that gains in  accuracy may require an exponential increase in computational effort and therefore quickly become unrealistic. Thus, in most problems we cannot hope to be able to write down an explicit solution of the mathematical model, and even a very accurate approximation of a solution may not be computationally feasible.

\section{What are the objects?}
As described at length, the objects of Newton's theory are modelled as
points that possess masses, but no extension. In solid state physics,
physicists are working with objects that have more internal properties
and structure and therefore might appear more
realistic. Electromagnetic and other fields are  objects of
physical theories with less material substrate. The nature of
subatomic particles is not so clear, even though one might grant them
some reality because they produce discernible effects \cite{Hac2,Fal}. When we go to
string theory, particles are nothing but excitation modes of
completely virtual ``objects'', the strings (see e.g. \cite{GSW,Pol,J5,J4,Zwie}). Such an excitation mode
is similar to a Fourier coefficient. There is nothing material about
it, and it is a good question in which sense that should be considered
as a physical object. Clearly, however, it is a theoretical object,
even though it might not possess any independent reality. 

Many biological structures are modelled as discrete entities, even
though they might be composed of atoms, molecules, cells, or
individuals in a complex manner. Mendel's laws are a good
example. Mendel conceived of genes as abstract entities, without any
clue about their physical implementation or realization. Although
Mendel's laws are not exactly valid, for many purposes this does not
matter. This leads to the question how such discrete entities can
emerge from some continuum at a lower scale. Conversely, large
ensembles of such discrete objects are often best modelled by
continuum models, as in population genetics, see e.g. \cite{HJT}. 

In any case, even though it has been discovered in molecular biology
that the material substrates of genes are nucleotide sequences, the
modern biological conception of a gene is more abstract than such a
nucleotide sequence. We refer to  \cite{SJ2} and the discussion in the journal
Theory in Biosciences about this issue \cite{SJ3}. Similarly, the other fundamental biological
concept, the species, is more abstract than the collection of
individuals composing a population (see for instance \cite{BrJ}). Thus, ``objects'' in modern biology are abstracta like genes or species that do not directly correspond to physical entities.

This is not fundamentally different in the social sciences, even
though there are fewer examples of successful mathematical models. At
present, the validity of the mathematical models of macroeconomic
theory is somewhat controversial. Likewise, it is an issue of debate
to what extent game theory can adequately capture human behavior (see for instance the discussion in \cite{KSJ}). On
the other hand, in recent years there has been definite progess in
modelling traffic dynamics with quasi-physical theories inspired by
microscopic particle physics, 
gas kinetics or 
fluid dynamics (see \cite{Hel}), notwithstanding the fact that individual drivers can
behave very differently from each other.

Another issue that has been debated already 200 years ago by social
scientists and mathematicians is to what extent the law of large
numbers provides not only an empirical, but also a conceptual basis
for mathematical modelling of social phenomena and dynamics, see \cite{Hac1}. When the
samples are large enough, statistical laws can lead to arbitrarily
precise predictions. This is the foundation of such domains as
demography or demoscopy. When we recall Boltzmann's description of
statistical ensembles or the above scale transitions, this appears no
longer fundamentally different from what has been discussed above for
physics and biology.

\section{Scales and levels}\label{levels}
In the preceding, I have already discussed the fact that the difficulty and feasabilitity of a description or a model of a complex system can be rather different at different levels. In this section, I want to analyze this issue more systematically, drawing upon the collaboration with Nihat Ay, Nils Bertschinger, Robin Lamarche-Perrin, Eckehard Olbrich and Oliver Pfante within the EU Project MatheMACS (Mathematics of Multilevel Anticipatory Complex Systems). Let me start with the terminological distinction between {\it levels} and {\it scales}. According to the definitions that we propose
\begin{itemize}
\item Scales refer to observables. Scales are determined by the measurement process, the data and their representations at different resolutions. For instance, concerning the dimension of length, we might consider the scales of micrometers, millimeters, meters, and kilometers, or other, smaller or larger ones.\footnote{In different fields, it may be different what is considered as a larger or a smaller scale. In geography, for instance, larger scale means higher resolution, that is, a smaller reduction factor. In other areas, like physics, a larger scale means the opposite. We shall follow the latter terminology. Thus, at a larger scale, many details from a smaller scale may disappear whereas larger structures might become visible.}
\item Levels refer to descriptions or models. Levels arise in the modeling process through the identification of entities that lend themselves to useful analysis. Thus, we could model a biological system at the molecular, cellular, tissue or organismal level. 
\end{itemize}
Clearly, the two concepts are interdependent. Levels rely on
measurements taken at characteristic scales. They also rely on a
choice of observables that are measured. Scales can be chosen
relatively arbitrarily, as long as the relevant  measurements can be
performed, but a useful choice of scale should depend on the
identification of a level. Levels, however, should be distinguished by characteristic properties of the model they allow. That is, at a given level, particular regularities should arise that do not pertain  at other levels. 

First of all, there are the systems that do not possess any
characteristic scale. Such systems have been objects of research in
both mathematics and physics. One of the concepts proposed here is
that of self-similarity. A self-similar structure looks the same at
any scale. Such structures were popularized by Mandelbrot \cite{Man}
under the name of fractals. The relevant mathematical theory had
already been created earlier, by Hausdorff \cite{Haus} in 1918.  Hausdorff had developed a general concept of dimension,
and many fractals possess a non-integer Hausdorff dimension. In particular, fractals are in some sense highly irregular. Nevertheless, there also exist such fractals with an integer dimension, and so, mathematically, a more subtle mathematical characterization is necessary, see \cite{Stef} and the references therein. The concept of a fractal also comes up in the theory of chaotic dynamics where the attractors and the boundaries of their domains of attraction could have such a fractal structure. The corresponding discoveries about iterations of polynomial maps were also made in 1918/19, by Julia and Fatou, see \cite{Dev} for a more recent description. In the physics literature, scalefree structures are usually characterized by a power-law behavior, as opposed to an exponential decay of correlations. That is, we essentially find correlations at any scale. It turns out that such power-law behavior is rather ubiquitous, from the financial data analyzed by econophysicists to the degree sequences of empirical networks in different domains. From a more theoretical perspective, such systems can even serve as some kind of universal models, when arising as limits of renormalization group flows (see e.g. \cite{Cardy}), and they are basic constituents of conformal field theory and the theory of critical phenomena, see \cite{DMS,ZJ}. 

Of course, no real system can be self-similar at all scales or possess
correlations of all orders. Such properties can only hold within a
certain range of scales. Nevertheless, as asymptotic idealizations,
 systems with such properties can be quite useful as theoretical
 tools.\footnote{Batterman \cite{Bat} and Lenhard\cite{Len} emphasize the fact that such idealized
 systems that arise as asymptotic limits in some theory are employed
 as models for empirical systems.} Scalefreeness, however, is not an
aspect I want to focus upon here. Let me rather return to a system
with well-defined objects at some particular scales. For instance,
take a couple of Newtonian particles, at some scale, for instance gas
molecules or celestial bodies. When there are only a few of them,
their kinetic or dynamical interactions can be described by a
mechanical model. For the celestial bodies, that is, when gravity is
the active physical force, we have already seen above that for more
than two particles, an explicit solution of the equations of motion
in general cannot be achieved. But even worse, when the numbers of
particles becomes large, like the stars in a galaxy or the gas
molecules in a container, the description in terms of classical
mechanics is no longer feasible, and we need to turn to a less
explicit statistical description \`a la Boltzmann. Here, a particular
configuration of, say, gas molecules, given in terms of their
positions and momenta, is a microstate, but at the macrolevel, one
only has collective observables like entropy or temperature. These
observables are not meaningful at the level of the individual
particles, because they represent statistical averages. More
precisely, any macrostate can be the result of many different
microstates, and the more microstates underlie a given macrostate, the
higher the latter's entropy, that is, the more likely it is to
occur. (There are some subtle issues here concerning the
interpretation of the probabilities involved, see \cite{Jay}, but we
do not enter into those here.) In summary, here a deterministic
description at the level of individual particles becomes unfeasible
and useless due to their large numbers and the impossibility of
accurately measuring all their positions and momenta, and instead a
statistical description at a higher level is used. It can also be the
other way around, that a statistical description at a lower level
yields to a deterministic description at a higher level in terms of
statistical averages. As fluctuations may average out, these averages
themselves may well obey deterministic laws. We see this, of course,
in the transition from the quantum level to that of micro- or
mesoscopic physics, but for instance also in cells when we go from the
molecular to the cellular level. Recall our discussion of the
Hodgkin-Huxley equations in Section \ref{physics}. 

Besides the dichotomy between deterministic and stochastic models,
there also is the dichotomy between discrete and continuous
models. When we move up or down the scales, or better, from one level
to the next, again the relation can go either way. Discrete particles
(atoms, molecules, stars, cells, individuals, ...) can give rise to
continuum models at a higher scale. However, out of underlying
continua, also discrete structures can emerge. Partly, these are
simply physical processes, for instance when interstellar dust
condenses into a star by the force of gravity. Partly, the discrete
entities are crucial constructs of models. The gene as a concept
in molecular or evolutionary biology emerges from some underlying
level of molecular interactions or reproductions of individuals in a
population, see \cite{JS}. The continuous dynamics of the
Hodgkin-Huxley model gives rise to the discrete event of a spike. Or
putting it somewhat differently, the neuron as a biological system
transforms an analogous input into a binary output, that is, it
chooses between the alternatives of spike vs. no spike. Thus, models
of information transmission in neural systems can operate with
discrete events as the carriers of discrete bits of information. From
a formal point of view, such phenomena are genuinely nonlinear, and
there exist theoretical concepts, like those of  critical threshold,
bifurcation, or phase transition, to analyze them, see for instance
\cite{J2,J6}.  Here, however, rather than  analyzing the dynamical
mechanisms that give rise to such discrete events in a continuous
setting, I want to utilize this as an example of a higher level
description in discrete terms with an underlying continuous dynamics
at a lower level. Interestingly, at the level of discrete events, the
spiking pattern of a neuron is often modelled as a Poisson process
(see for instance \cite{DA,GK,J7}),
that is, as a stochastic process, instead of a deterministic
one. Here, the reason why one switches from a deterministic to a
stochastic description is somewhat different from the case of the
Boltzmann gas. In the latter case, a statistical description is called
for because of the larger number of particles involved and the
infeasibility or impossibility to measure all of them. Here, in
contrast, we have a single event, the generation of a spike by a
neuron that at one level is modelled by a deterministic dynamical
system.\footnote{Interestingly, stochastic perturbations of this deterministic system produce genuinely nonlinear effects, see \cite{GTJ,TJ}.}  Also, this system does not exhibit chaotic behavior that would
make long term predictions impossible. Of course, near the critical threshold
above which the neuron spikes and below which it returns to rest, the
dynamical behavior naturally is unstable in the sense that arbitrarily
small perturbations or fluctuations can induce a switch from one state
to the other. At
another, more abstract, level it is instead modelled as a stochastic
process. This is not only simpler, because it dispenses us of having
to deal with the subtle nonlinear behavior of the Hodgkin-Huxley
system, but it also offers the advantage that we can now look at the
relation between the strength of the neuron's input and the single
parameter that characterizes the stochastic process, the firing rate
of the neuron. Put simply, the hypothesis would be that the neuron
fires more frequently on average when it receives a stronger input. (Whether, or perhaps more precisely, to what extent neuronal systems really employ a rate coding as
the preceding might suggest is an unresolved and intensely debated
issue in the neurosciences, but here, I do not
 enter that discussion, and rather refer to \cite{GR,WVL} for the information theoretical  analysis of neuronal spike trains.)

I now turn to the formal analysis of the relation between different
levels, following \cite{PBOAJ1}. (Let me also mention earlier results in \cite{ShaM,GNJ} and many other papers, and the case study
in \cite{PBOAJ2}.)

We consider a process
\bel{level1}
\phi:X \to X'.
\qe
This could, for instance, be the transition from the state of a system at time $t$
to its state at time $t+1$. We assume that this is a Markov process,
in the sense that knowledge of the state at prior times $t-1, t-2,
\dots $ does not contain any information beyond that contained in the
state at time $t$ relevant for the state at
time $t+1$. That is, the future is conditionally independent of the
past given the present. We also have an operator
\bel{level2}
\pi:X \to \hat{X}
\qe
that is considered as a projection, coarse graining, averaging, or
lumping. This simply means that the $\ \hat{}\ $ indicates a
description at a higher level, with less detail and resolution. The
question then is whether we can transform the transition directly at
the higher level, that is, whether there exists a 
 process
\bel{level3}
\psi:\hat{X} \to \hat{X}'.
\qe
that is self-contained in the sense that it does not depend on
anything in $X$ that is not already present in the upper level
$\hat{X}$. In other words, we ask whether we can close up the following
diagram (make it commutative in mathematical terminology \cite{J8}).\footnote{A perhaps somewhat technical point concerning the representation by this diagram: Often, one thinks of a projection as going down, instead of up, and one would then represent $X$ in the top and $\hat{X}$ in the bottom row. Since, however, we think of $\hat{X}$ as a higher, more abstract, level, we rather represent that higher level in the top row of our diagram.}

\bel{cat1-11a}
\begindc{\commdiag}
  \obj(1,3)[1]{$\hat{X}$}
  \obj(4,3)[2]{$\hat{X}'$}
  \obj(1,1)[3]{$X$}
\obj(4,1)[4]{$X'$}
 \mor{1}{2}{$\psi$}[\atleft,\solidarrow]%\dasharrow produces error message
  \mor{3}{1}{$\pi$}[\atright,\solidarrow]
\mor{4}{2}{$\pi$}[\atright,\solidarrow]
  \mor{3}{4}{$\phi$}[\atright,\solidarrow]
  \enddc 
\qe 
Expressed as a formula, commutativity means that
\bel{cat1-11b} 
\psi (\pi(x)) = \pi(\phi(x)) 
\qe
for all $x\in X$.

There exist several criteria to make this precise. 
\begin{enumerate}
\item[I] \textbf{Informational closure}: All we need to know to
  determine the state $\hat{X}'$ or to predict its statistical
  properties is already contained in the state
  $\hat{X}$. The higher process is informationally closed, i.e. there is no
information flow from the lower to the higher level. Knowledge of 
the microstate will not improve predictions of the macrostate. The upper level process is self-consistent in the sense that it does not need
to perpetually draw information from or about the lower-level states.
\item[II] \textbf{Observational commutativity}: It  makes no difference
  whether we perform the aggregation first, 
and then observe
the upper process, or we observe the process on the microstate level, and then lump together the states. In that sense, the upper level process seems  autonomous, as once initialized, it appears to unfold on its own, without needing any further updating by details from the lower level process. 
\item[III] \textbf{Commutativity} in the sense of \rf{cat1-11b}: There exists a transition kernel $\psi$ such that the diagram  \rf{cat1-11a}
 commutes.
\item[IV] \textbf{Markovianity}: $\hat{X}, \hat{X}'$ forms again a
  Markov process. We recall that we assume that $X, X'$ yield  a Markov process, that is, the current state of lower level process contains everything needed to compute its next state, and it does not need to draw upon past states any further, as these have transmitted all relevant information to the current state. But  it does not follow in general that the upper process is Markovian as well, and so, IV indeed is a nontrivial condition. It could compensate a lack of information about the current lower level state by structural constraints in order to also utilize information from past states.  In particular, at the upper level, information about past states can improve its prediction about the next state, whereas such information is not useful at the lower level.  See the example in \cite{PBOAJ2}, and also the discussion in Section \ref{language}. 
%\item[V] \textbf{Predictive efficiency}: An emergent  level corresponds an efficiently 
%predictable process. 
\end{enumerate}

We show in \cite{PBOAJ1} that I implies both II and IV, and II
implies III, whereas IV does not imply III in general. Also, for a 
deterministic process, but not necessarily for a stochastic one,
conversely III implies II which in turn implies I. 
III of course is a formalization of I and II, and as we show, in the deterministic case, they are equivalent. In the stochastic case, the notions turn out to be somewhat different, as information concerns statistical properties, and these are different from the results of particular observations. IV, of course, is a condition that is meaningful only for stochastic dynamics.

\section{The structures of data sets}
Let us look at an example that may seem very simple at a first glance,
a collection of minerals in a museum. How should the specimen be
arranged and ordered? There are many possibilities:
\begin{itemize}
\item according to geographical origin, that is, where they have been
  found
\item by physical properties, like color, hardness, 
  size, 
  etc.
\item by chemical composition
\item by geological age
\item by their process of generation, like volcanism, pressure,
  fossilization, etc. 
\end{itemize}
Thus, with a heterogeneous data set, a good ordering scheme may be
neither easy to find nor unambiguous. In the case of the minerals, however, a solution came from an important discovery in physics. Max von Laue, on the basis of his work with Walter Friedrich and Paul Knipping on X-ray diffraction in crystals,  detected that crystallographic lattice structures could be determined by X-ray spectroscopy (see for instance the account in \cite{vL}). And this implied that minerals can be classified in terms of the symmetries of their crystallographic structures. In fact, this  classification by structural symmetries is the same principle as that used by 19th century biologist Ernst Haeckel \cite{Hae} in a very different realm, for the  classification of radiolariae, protozoa that provide much of the plankton in the oceans and that produce mineral skeletons  of highly symmetric  geometric shapes. This is another instantiation of the fact that one and the same mathematical concept can be applied in very different data sets. 

A look at the history of science may also be illuminating here. In
particular during the 17th century, there was a widespread belief in a
natural order (as imposed by God) within which each individual object
or class of objects would find its place. Perhaps the most prominent
application of this principle was to the classification of plant and
animal species. That is, there was an abstract scheme where each
species found its place in, and conversely, every position within this
scheme had to be filled by precisely one species, already known or
still to be discovered. The important step of the botanist Carl von
Linn\'e then was to infer the natural order from a careful
inspection and systematic comparison of all known plant species,
rather than from some a priori reasoning. Still, the difference to the
Darwinian concept of evolution is profound. According to Darwin, there
is no natural and systematic order, but rather the collection of
species is the result of a partly contingent historical process. An
interesting reference is \cite{Amun}.

\section{Data without underlying objects?}

When we start with data, then an object behind these data is first a hypothesis, a construct. A well discussed example are the gigantic data sets created by modern particle accelerators  where the elementary particles that are supposed to produce those data are postulated by theory, but where the only evidence for them has to be extracted from those data (for a conceptual discussion, see \cite{Fal}). 
More generally, physics studies phenomena and data, like the scattering matrix in quantum mechanics that records the input and output of a quantum mechanical system. 
Physicists do not think in terms of objects, in contrast to the picture
of Newtonian physics depicted above. In fact, precisely these
scientists are those that drive the data analysis approach also in other domains. For
instance, neoclassical economic theory has developed elaborate models
for describing the functioning of an (idealized) economy, but then
came the so-called econophysicists that simply took the data from
financial and other markets and searched for statistical
regularities \cite{MS}. They used methods from nonlinear time series analysis,
for instance. Those methods are not sensitive to the origin of the
data, nor do they depend on any object oriented models. Rather, they
look for intrinsic regularities in possibly chaotic dynamics
\cite{KS}. For instance, they identify (chaotic) attractors and
determine their (possibly fractal) dimension. These, however, are not
ad hoc methods; they are based on deep insights from the theory of
dynamical systems. For example, the center manifold principle, see
e.g. \cite{J2}, (which, for instance, is the basis of Haken's slaving
principle \cite{Hak}) tells us that in a dynamical system, typically
most directions are quickly relaxing to their equilibrium, and the
essential aspects of the dynamics are determined by very few slow
variables.

In fact, it had already argued by Koopmans \cite{Koop} in his review
of a book by Burns and Mitchell \cite{Burns} that  measurements of
economic quantities without underlying economic concepts  are
useless; and that time, the criticism was directed against a use of
econometric techniques without a theoretical framework grounded in
economics, but the same type of criticism would also, and perhaps even
more forcefully, apply to the approach taken by econophysics.

We should also remember that Newton's theory that can derive and
predict the motion of celestial objects is a rather exceptional feat
even in the history of astronomy. Astronomy started with data
collection, and these data then were typically organized by models
with an ontological status that was dubious at best. The old
Babylonians simply interpolated between observational data by assuming
that the positions of the celestial bodies varied linearly between the
extrema and then suddenly changed  direction, see the systematic studies of Neugebauer \cite{Neu1,Neu2}, who speaks of zigzag functions to describe this interpolation scheme. The Ptolemaic
model as described in the Almagest employed more and more epicycles without really accounting for
their ontological status. The Aristotelian system of the celestial
spheres was more explicit in that regard, but more difficult to
reconcile with the astronomical data. (In fact, Ptolemy also not only tried to represent data as in the Almagest, but also developed a planetary hypothesis in which Aristotle's celestial spheres were fattened to provide space for the epicycles needed in Ptolemy's system to account for the astronomical data.) And Kepler had to labor
painstakingly for many
years through the data compiled by Brahe. He first tried a 
preconceived scheme,  an organization of the solar system in
terms of Platonic solids, which he then  abandoned because it did not fit the data. In particular, he was concerned with the orbit of Mars. He tried many variants before he finally arrived at fitting the orbit as an ellipse.  Thus, by  trying to account for the data, in the end  he succeeded in discovering his laws of
planetary motion. As already discussed, the laws that Kepler had
empirically discovered were then derived
by Newton from his law of gravity, that is, within an explicit
physical model. Of course, without knowing Kepler's laws, Newton might
not have found his theory of gravity.

I should concede, however, that this picture of large data sets that
are first collected without an adequate theory and only subsequently
inspire deep physical theories is not without exceptions. At least one
exception springs to mind, Einstein's theory of general
relativity. This theory was developed on the basis of an abstract
principle, general covariance, and was only subsequently tested
against and confirmed by empirical data. 

Also, the traditional view of physics is that it does not collect more
or less 
arbitrary data, but conducts specific experiments whose results
acquire their meaning within an established theoretical framework. 

On the other hand, however, we should also discuss the approach of
Alexander von Humboldt. The novelty of his expeditions rests in the
systematic and comprehensive collection of all kinds of data with all
the measuring devices available at his time, at a time when there was
still very little theoretical understanding of the processes shaping
the ecology of the earth. His expedition in South America took some years, but the
evaluation of the data collected during that expedition took him
several decades. His work then launched the scientific discipline of
physical geography and helped and inspired several other scientific
fields, even though no single coherent and encompassing theory emerged
from his data. Thus, Humboldt had gathered huge amounts of data, but these data did reveal only very few coherent patterns, like the dependence of climatic zones on altitude and  latitude. In other words, not only was the model missing, but also the data analysis largely failed in discovering general structures.

\section{Towards an abstract theory}
In fact, it is a quite general finding that there exist profound
analogies between models in physically very different domains. It is
actually one of the main driving forces of mathematics to consider the
corresponding structures and relations abstractly and independently of
any particular instantiation and to work out general theories. This
will then make it possible, in turn, to apply these mathematical
theories to new domains. It is then irrelevant whether such a domain
constitutes some independent reality or whether one simply has a
collection of data. For instance, similar statistical phenomena show
up in quantum mechanics, in the analysis of biological high-throughput
data, or in the description of social phenomena. 

Understanding the laws and regularities of the transition between
scales (see e.g. \cite{PS} for the mathematical background), or in a more ambitious formulation, a mathematical approach
towards the issue of emergence (for instance, \cite{JBO}), is a theoretical challenge that
transcends individual domains. Nevertheless, it should lead to
fundamental insight into the structure of reality,  at an abstract level. 

In particular, it is a fundamental question in the theory of complex
systems to what extent general laws apply across different domains and
disciplines, and where the applicability of a general theory ends and
a more concrete modelling of the details of the specific system is
required. For instance, which analytical concepts and mathematical
tools apply simultaneously to cells, neural systems, psychological
systems, and societies, or
at least to several of them, and
where do the specific peculiarities of each of these fields enter? 

In a different direction, we may ask for a theory of structure in
high-dimensional spaces. As will be explained below, even though the
data may possess a large number of degrees of freedom, these degrees
of freedom typically do not vary completely independently in large
data sets, but rather obey some nonlinear constraints. Thus, we are
dealing with intrinsically lower dimensional geometric structures in
high dimensional spaces. Since the particular structure will usually
not be known before having analyzed the data, we need to consider
spaces of such structures. This leads to new challenges for the
mathematical field of geometry. These spaces will then carry a probability
density, telling us about the a-priori likelihood of finding them in a
particular data set. In statistical terms, we then have some prior
hypothesis about the data set, and such a prior then has to be turned
into a posterior by Bayes' rule on the basis of the data observed. Of
course, this then should be done in an iterative manner. In a
different direction, we may employ methods from algebraic topology or metric geometry in
order to discover robust qualitative features of specific data sets,
see \cite{Car,BHJKS}.

\section{The challenge of big data}\label{big}
Both the abundance of huge data sets in almost all disciplines and the
availability of the computational power to formally analyze them are
quite recent developments. We have the high-throughput data in
molecular and cellular biology generated by gene sequencing,
microarrays or various spectroscopic techniques, the imaging data in
the neurosciences, the email or movement data of mobile phone users,
the data of transactions in financial markets at millisecond
resolution, the linking patterns of the world wide web, and so on. 

And almost uniformly across disciplines, we see a major transition
from model driven to data driven approaches. The latter approaches
typically depend on the statistics of large data sets. These
statistics are utilized  automatically, and statistical regularities
within large data sets are only used in an implicit fashion, without
making them explicit. Examples abound. In molecular biology, the
problem of protein folding had for a long time been approached by
explicit physical models. The problem consists in predicting the
three-dimensional folding pattern from the linear sequence of the
amino acids constituting a polypeptide on the basis of molecular
attractions and repulsions. As the energy landscape is quite
complicated, there are typically many metastable states, and finding
the configuration of minimal energy, i.e., the configuration
supposedly assumed by the polypeptide in the cell, therefore is
computationally quite difficult, even though powerful Monte Carlo type
schemes have been developed. Recently, however, it turns out that the
best predictions are achieved by data bank searches without any
physical model in the background. One simply compares the sequence at
hand with those sequences where the folding pattern is already known,
for instance by X-ray cristallography, and then makes a statistical prediction
based on sequence similarities.  

Perhaps the preceding becomes clearer when we consider a hypothetical approach to weather forecast. A model driven approach develops a detailed dynamical model of cloud formation and movement, ocean currents, droplet formation in clouds, and so on, and measures the many parameters of such a model, or perhaps also tries to fit some of them  on the basis of observations. Such a model would be given by a system of coupled partial differential equations (PDEs), for which one has to solve an initial value problem. The initial values are determined by current measurements on a  grid of measurement stations that is dense as possible, naturally with a higher density  on land than on the oceans. Since the dynamics described by the PDE model tends to have chaotic aspects (see also the discussion in Section \ref{valcomp}), the precision of the measurements of the initial values is of utmost importance in order to have somewhat accurate predictions for a few days. Likewise, the details of the model and its parameters are crucial. In contrast, a data driven approach would like depend on accurate measurements, but it would then try to identify those constellations in the past whose values are closest to those presently recorded, and then use the known weather dynamics from the past for those constellations to predict the weather derived from the current values. This would simply require large data bases of the past weather recordings, and perhaps some scheme of taking weighted averages over similar constellations in the past, but it would not need any model of weather dynamics. Such an approach should then naturally be expected to improve as  the data base grows over time. 

In the geosciences, the standard scheme consists now in estimating 2-point correlators  from noise correlations at different locations (see e.g. \cite{RLS,GP}), instead of using physical models of, for instance, wave propagation in geological media. The noise sources can be physical or caused by human activities. While the differences matter at a technical level because of different noise characteristics, this does not affect the principle of the method. 

Or to present another example that
plays a prominent role in \cite{Nap}, microarrays record the
simultaneous expression of many different genes in a particular cell
condition by exposing the corresponding RNAs in the cell
simultaneously to an array
of pieces of complementary DNA sequences and then simply recording
which of those DNA sequence pieces find RNA partners. The underlying
biological rationale is the following. The DNA pieces are part of the
genome that is selectively transcribed into RNA which may then be
further translated into polypeptides, the building blocks of
proteins. Proteins carry out most of the essential operations in a
cell. Therefore, mechanisms of gene regulation should ensure that
precisely those proteins are manufactured that are needed by the cell
in a particular situation. Therefore, precisely those pieces of DNA
should be transcribed into RNA that encode the right proteins (there
are some problems with this oversimplified account, see \cite{SJ2},
but the microarray technology happily ignores them). Thus, a
microarray tests the expression patterns and intensities of many DNA
segments -- which stand for genes in this simplified model --
simultaneously, because the transcribed RNAs bind to the complementary
DNA pieces offered by the microarray.

In computer linguistics, for the purposes of automatic translation,
models grounded in syntactic and semantic theory are replaced by
statistical techniques that simply utilize cooccurence patterns of
words in large corpora. In particular, Google uses so-called $n$-gram models for the purpose of automatic translation. This simply means that one derives the relative frequencies of strings of $n$ words from databases containing trillions of entries. $n=5$ is a typical value. That is, the meaning of a word -- insofar as
one should still speak about meaning here, in the absence of any
semantic concepts -- is determined
by the environment, that is, a few words preceding and following it, in which
it occurs. We shall analyze the conceptual shifts that this implies in more detail in Section \ref{language}.

In economics, as already mentioned, techniques from statistical data
analysis, like nonlinear time series analysis, are applied to
financial data sets in order to find subtle patterns that do not
follow from theoretical models. This is the new field called
econophysics. Recently, econophysicists also started to develop
economic models. They thereby depart from a purely statistical approach that solely analyzes economic or financial data and now build models of economies or financial markets themselves. Their models, however,  are in contrast to the classical
economic models that employ a
so-called representative agent. The latter stands for  models with many identical typical
agents that are ideally analytically tractable. The agent based models
of econophysicists instead utilize very simple, but possible diverse
and heterogeneous agents that fit well to large scale computer
simulations. That is, instead of a single type of agent that might be rather
sophisticated and grounded in economic theory, here many simple agents
are employed. These agents depend on a couple of parameters, and the
parameter values can and typically do differ between the agents. The representative agents of economics are usually assumed to be fully rational -- the argument being that non-optimal agents are quickly exploited by their more clever competitors and thereby driven out of the market. The diverse agents of agent based models are not at all assumed to be rational. They rather follow relatively simple empirical rules, and the purpose of the models is to uncover the collective effects of such behavior of many agents through systematic computer simulations. The
validity of the simulation results usually remains unclear,
however. An important issue is that the models of (neo)classical economic theory strive for the ideal of exactly solvable equations, or at least for describing the economy by a small and carefully specified set of explicit equations. On the basis of these equations and their solutions, they want to achieve analytical insights into the working of the economy. The agent based models of econophysicists, in contrast, employ much less rigid models, possibly  with many parameters, and large numbers of equations that can only be numerically solved. Analytical insight no longer is the foremost aim, and one rather wants to identify unexpected nonlinear effects and critical transitions between different regimes triggered by small variations of crucial parameters.  

Even in elementary particle physics, after the confirmation of the
Higgs boson, that is, a prediction made by the standard model, at the
LHC, in the future one will probably move towards the automatic
analysis of large scale scattering data in order to find unpredicted
events that may not fit into any current theoretical model. In any
case, already for finding evidence for the Higgs boson, to a large
extent automatic data analysis methods have been employed. Such
methods may find patterns of correlations or untypical events in huge
data sets by completely implicit methods. As described and analyzed
for instance in \cite{Fal}, the empirical evidence is gathered in four
steps. At the basis, there are position measurements. Adjacent
position measurements are then combined into tracks. Adjacent tracks
in turn are combined into events. Finally, statistical ensembles of
scattering events contain resonances. All this is done in a completely
automated way. The tools have no special affinity to particle physics,
even though particle physicists are among the people making the most
advanced use of them. For instance, so-called
neural networks  encode the patterns of their training sets in a
rather indirect and 
completely implicit manner in synaptic connection weights (see e.g. \cite{Bi}). When they
are subsequently applied to the real data, they produce corresponding
associations which may then be interpreted as patterns in the data
set.  Neural networks and other such schemes find applications in a
wide range of data domains. In fact, the methods themselves also
change, and for instance, neural networks are sometimes replaced
by other methods like support vector machines.\footnote{Support vector machines are efficient classifiers that use a high-dimensional linear feature space (see  \cite{CST,SS,SC}).}

This phenomenon, the transition from model to data driven approaches,
occurs not only in individual and specific domains, but also at the
formal level. The field of scientific computing currently undergoes a
similar transition. The field of machine learning is concerned with
the development of techniques for the automatic analysis of large data
sets. Statistics, the science of finding structure in data,  is moving into a similar direction, and in fact a
profound convergence between machine learning and statistics seems to
take place. In particular, the computer intensive Bayesian methods
take over much of more traditional parametric statistics. Often, such
methods are combined with stochastic search algorithms like Monte
Carlo methods or Markov type models. Lenhard\cite{Len} also points out
the interdependence between data dynamics and computational
modelling.

Big data are said to be characterized by large or even huge values of the three {\bf V}s, that is,
  \begin{itemize}
  \item {\bf V}olume: often petabytes ($1000^5 =10^{15}$ bytes)/day, possibly more
\item {\bf V}ariety: heterogeneity of data types, representation, and semantic interpretation
\item {\bf V}elocity: arrival rate and reaction time 
  \end{itemize}
Processing big data requires adapted strategies and methods, and it can be decomposed into five phases
  \begin{enumerate}
  \item Acquisition and recording: filtering, compression, metadata generation
\item Information extraction, cleaning, and annotation
\item Integration, aggregation, and representation; data base design
\item Analysis and modeling; querying and mining the data
\item Interpretation; visualization; possibilities of interaction of human observers with machine processing
  \end{enumerate}

For our purposes, phase 4 is the most relevant. 
Although  computer power is rapidly growing, and cloud
computing or even access to supercomputers becomes ever more
available, large data sets still may give rise to the ``curse of
dimensionality''. This means that computer time will increase
exponentially with the number of degrees of freedom and therefore
quickly exceed even the capacities of supercomputers, unless clever use
of specific structures within a data set is made. Therefore,
corresponding data analysis techniques are being developed. An example
is compressed sensing, see \cite{Do,CRT,FR}. Here, the idea is to use a specific sparsity
assumption, for instance that a large and complicated acoustic data
set might be generated by a small number of sound sources
only. For instance, the sound sources might be humans carrying on conversations in a crowd. Formally, sparsity means that a high-dimensional vector or
matrix possesses only few entries that are different from 0. One does
not know beforehand, however, which are the nontrivial ones and how
many of them there really are.  Thus,
there is a not very explicit, but still rather powerful and
constraining structural
assumption. Or, in a similar vein, one might assume that the data
points one has in some high-dimensional space are in fact constrained
to some low-dimensional manifold (\cite{BN}). This low-dimensional manifold,
however, need not be a linear space, and therefore, even though it
possesses only few intrinsic degrees of freedom, it may still stretch
into many of the ambient dimensions. Thus, the data are supposed to
possess some rather constraining type of regularity, but again, this
regularity 
is not explicit, and it is the task of a machine learning scheme to
use such an abstract assumption in an efficient manner. When
one incorporates such an assumption, one may drastically reduce the
number of computations needed to identify the underlying structures
from the data set, for instance the sound sources in the above
example. The mathematical techniques required are typically nonlinear
and therefore on the one hand more flexible in discovering and utilizing
hidden regularities and patterns in the data, but then, on the other
hand, require more specific, but ideally still automatic, adaptations
to the data set at hand. 

In many respects, human cognition is still superior to automatic data
analysis techniques. In order to utilize the power of human cognition,
however, it is necessary to convert the data into a format that is
familiar to humans. This leads us into the domain of visualization, a
newly emerging scientific discipline between computer science,
mathematics, and psychology. Thus, the aim is to develop formal
methods that can convert a data set into a form in which humans can
easily discern patterns. At a more advanced level, the combination of
machine learning tools and human perception may become interactive,
which 
was the goal of the EU funded research project CEEDs (The Collective Experience of Empathic Data Systems).

\section{Big data and automatic translation, or a paradigm shift from
  linguistic theory to language processing}\label{language}

Modern linguistic theory and philosophy is founded upon basic
oppositions, for instance between
\begin{itemize}
\item {\it langue} (the abstract system of a language) vs. {\it
    parole} (the concrete utterance) (de Saussure \cite{Sau})
\item {\it diachronous} (across time) vs. {\it synchronous}
  (simultaneous) (de Saussure \cite{Sau})
\item {\it competence} (the ability for the correct syntax of one's
  native language) vs. {\it performance} (the actual production of
  utterances) (Chomsky \cite{C1})
\item {\it deep structure} vs. {\it surface structure}, i.e., a
  language independent representation in an abstract structure is
  transformed into a sentence according to  the syntax of a specific
  language (Chomsky \cite{C1}), and this transformation obeys the
  general rules of abstract grammar; the latter have become more
  general and abstract themselves in the course of Chomsky's work
  \cite{C2,C3}
\item {\it spoken} vs {\it written} language, again from de Saussure
  \cite{Sau} and emphasized more recently by Derrida \cite{De}. 
\end{itemize}
Perhaps the most important opposition is that between de Saussure's
\cite{Sau}
\begin{itemize}
\item {\it paradigmatic} alternatives and {\it syntagmatic} series;
  this means that a word at a given position in a sentence can be
  chosen paradigmatically from a list of words that could
  grammatically occupy this position, but has to obey the syntactic
  rules of the sentence. This is also an opposition between {\it
    absence} and {\it presence}, or between {\it selection} and {\it constraint}; the alternatives that have not been
  chosen are absent, but the syntactic constraints are present through
  the other words in the sentence. 
\end{itemize}

These oppositions are (more or less) formal, but not directly
mathematical. They gave rise to a formal theory of language, and for
some time, one attempted to utilize that theory for purposes of
automatic translation. The idea was essentially to automatically
infer the grammatical relationships within a given sentence, 
that is, the dependencies between the different words and grammatical
tokens and the internal references, like the referents of pronouns and
anaphora. That grammatical structure could then be transformed into
the corresponding structure of another language, and the meanings of
the individual words could be correlated with the help of good
lexica. For instance,  for some time \cite{PS} was popular as a suitable grammatical
theory for such purposes. All such attempts, however, have more
recently been brushed aside by the big data approach to automatic
translation, as developed and pushed in particular by Google. That
approach completely ignores, and often even ridicules, linguistic
theory, and rather draws upon correlations in huge linguistic
corpora. My purpose here is to analyze the underlying conceptual shift
and to describe what mathematical structures are behind this
approach. Those mathematical structures are completely different from
those developed from the context of formal linguistics. And the people
working on automatic translation within this paradigm express little
interest, if at all, in the conceptual basis of their
endeavor. Nevertheless, that  will be
important for  this essay.

Thus, the oppositions sketched above no longer play a role. Instead,
we see the {\it corpus}, that is, a data base of texts in the
language(s) in question, as the single basic element. The corpus
scales by size. Thus, a small corpus might be seen as corresponding to
 {\it parole}, whereas big corpora can approach the {\it langue} side of de
Saussure's dichotomy, and the same applies to Chomsky's dichotomy
between {\it competence} and {\it performance}. More precisely, there
no longer is any such abstract thing as {\it competence}. The possibly
huge collections of linguistic data are all there is. Corpora do not care much whether they have been 
assembled from contempory, hence essentially simultaneous texts or
whether they result from scanning texts over longer periods of
time. The only relevant criteria are of a practical nature, digital
availability and computational capacity. Also, corpora can as well be based on automatically recorded
spoken language as on written texts. The difference is again largely
irrelevant, or at least plays no basic conceptual role. (More
precisely, what remains is the technical difference between off-line and
on-line processing and translation.)

More importantly and interestingly, the opposition between de
Saussure's {\it paradigmata} and {\it syntagmata} is also
resolved. This best explained through the mathematical concepts of
stochastic processes and Shannon information (see e.g. \cite{CT,SW,MK,J2}). De Saussure talked about
alternative words at a given position in a sentence. This is
qualitative, but not quantitative. Information
theory, in contrast, would quantify the probabilities of different
words to occur at such a position. These probabilities then will not
only depend on the abstract properties of that position, but also on
the concrete words before and behind it in the sequential string
forming the sentence or the text. That is, we do not just have
probabilities for the occurrences of words at a given position, but we
rather have transition probabilities from a word, or more precisely, a
segment of a few words, to the next. In fact, automatic translation
today mostly works with pentagrams, that is, overlapping strings of
five words.\footnote{In the terminology of statistical physics, the correlation length of word sequences in texts becomes relatively small after five words. In fact, one would expect that it never becomes exactly zero, that is, there do exist correlations of arbitrary length, whatever small. Thus, in technical terms, when moving from a word to subsequent ones in a text, we  have a stochastic  process that does not satisfy a Markov property for strings of words of any finite length.} Thus, we no longer have de Saussure's qualitative
opposition between selection and constraints, but both aspects are
combined and made quantitative within the concept of a stochastic
process. (In fact, this is a slight oversimplification. A stochastic
process would occur if one has to guess or produce the next word on
the basis of those preceding it. Thus, a process describes a temporal
sequence as in actual speech. A corpus, however, is not a process
unfolding in time, but is simultaneously given as a whole. Therefore,
the probabilities are determined not only by the preceding words, but
also by the following ones. The underlying mathematical principle,
however, is not fundamentally different.) Thus, while the big data
approach to automatic translation does not care about an underlying
conceptual structure, its analysis nevertheless sheds some light on
the limitations of formal language theory, and it points towards
mathematical structures that replace the qualitative oppositions
by quantitative probabilities. These probabilities are nothing but
relative frequencies of pentagrams, i.e., certain word constellations,
and they can be automatically computed from the corpora at hand. The
larger the corpus, the more accurately such probabilities can be
determined, in the sense that further additions to the corpus will
likely 
change them only very little. In a
certain sense, this approach asks more precise questions than formal
linguistics. When analyzing a temporal sequence, like a spoken text,
it would not simply ask at each instant ``what could come next?'', but
rather ``how well can I guess what comes next?'', and the latter is
quantified by Shannon's information \cite{SW}. 

Of course, this approach to automatic translation has its limitations. It cannot capture long range dependencies. In a nutshell, the basic assumption underlying the approach is that the transition between overlapping pentagrams satisfies a Markov property, that is, no information from more than five words back in the sequence is needed  for the probability of the next word. From the syntactic perspective, this for instance does not adequately capture the structure of the German language where the different components of the verb can be separated by many words in a sentence. From the semantic perspective, the meaning of some sentence in a text may refer to other, much earlier parts of that text, or even to a context outside the text itself. In formal terms, let us recall also our discussion about Markovianity at different levels in Condition IV in Section \ref{levels}. 

One of the central, but to a large extent still unresolved, issues of
linguistics, that between syntax and semantics, between structure and
meaning, is simply bypassed by automatic translation. A human
translator would first try to extract the meaning of a text and then
express that meaning as well as she can in the other
language. Automatic translation, in contrast, simply transforms one
structure into another one.

\section{The issue, again and hopefully clearer}
We are faced with the following alternative to which no general answer
can be given, but which needs to be evaluated in each individual
situation.
\begin{enumerate}
\item The significance of  data collected
  depends on the specific context from which they are taken, and they
  cannot be fully understood without that context.
\item Data sets typically possess internal structure, and much of such
  structure generalizes across different disciplines, domains and
  contexts, or is at least accessible to context independent methods. Identifying and understanding that structure with formal
  tools will then in turn enable us to use the data to learn something
  new and 
  insightful about their context. 
\end{enumerate}
This is nothing but the old alternative between specificity and
generality, as, for instance, already described by Kant. It is easy to deplore too much of an emphasis on either
side and support this by supposedly misguided scientific case studies, and to declare a
one-sided practice as agnostic or unscientific. The real challenge
consists in finding the appropriate balance between the two aspects in
each case. Currently, big data sets offer new opportunities for formal
methods and computational models, and this may shift the balance for a
while. 

Concerning the role of mathematics, this might transcend the
alternative between a mathematics of content that establishes and
analyzes relations between the concepts of a theory and a purely
auxiliary mathematics that is relegated to the  role of
data handling and preprocessing. What should emerge rather is a
mathematics that develops new abstract concepts for analyzing and
representing data spaces. For instance, a partial merging of the domains of
statistics, high dimensional geometry, information theory and 
machine learning might take place. Such a mathematics would help to
detect structures, as opposed to either formalizing structures
proposed and 
developed by other scientific disciplines or to offering tools for fitting given data
into such a structure.  This process can only enhance the role and the
importance of mathematics, as it becomes liberated from being
instrumentalized by conceptualizations that are externally imposed. Instead of accepting models from particular domains, mathematics would itself propose abstract metamodels, like sparsity, smoothness, or symmetry.

\section{Some consequences}
The preceding has important implications for the role of models and
hypotheses in the scientific process. From the perspective of
Bayesian statistics, we begin the process of scientific inquiry with
prior hypotheses constructed by us, in whatever way we may find
reasonable or on the basis of whatever prior experience we may
have. During the process, such a prior hypothesis gets transformed into
a posterior one on the basis of the observations made or the data
acquired. This is achieved by applying a fixed rule, that discovered
by Bayes. The role of models then is relegated to provide sets of
parameters that have to be fitted to the data. This is in stark
contrast to the ideal of a model in physics. Such a model should
contain only very few, and ideally no free parameters at all that are 
not theoretically determined but need to be measured. 
Even for such a nearly ideal model, one may question whether empirical
adequacy should constitute a proof of the correctness of the model. In
a model with many free parameters that need to be fitted to the data,
certainly empirical adequacy can no longer  count as a proof of the correctness of the
models, and indeed, in complex situations, the term ``correct model''
may lose its meaning entirely. {\bf There no longer is such a thing as a
correct model.} There may only be a superior fit of the parameters to
the data collected according to unknown probabilities. Alternatively,
we might be seeking regularities in data on the basis of certain
structural hypotheses, as described in Section \ref{big}. Again, such
structural 
hypotheses like sparsity do not come from the data at hand, but are
applied by us in order to get some handle on those data. In a certain
sense, they constitute a structural prior, although that prior is
usually not updated in a Bayesian manner. 

While all this may sound rather agnostic, in practice one might get
quite far with those schemes of Bayesian updates, parameter fitting,
and structural hypotheses. It is challenge for mathematics to analyze
this issue theoretically. \\

Putting it somewhat differently: Data analysis depends on prior
structural assumptions. That could be the prior of a Bayesian
approach, or it could be the choice of a family of models as in
parametric statistics. It could be a general structural hypothesis
like sparsity, a low-dimensional manifold (for instance arising as the
center manifold of a dynamical system), certain invariances or
symmetries etc. This is the epistemic side. At the ontological side,
is there anything underlying the data matching those assumptions? One
may argue that these assumptions are our constructions, and that
therefore there is no guarantee of any correspondance with the data
source. However, a simple theory can sometimes enable us to discover
specific phenomena that would otherwise not emerge from beneath the
heap of data. One may also argue, as I have tentatively done in this essay,
that in some situations, such assumptions could be justified by the
structures from which the data are derived, but there is no guarantee
for that. The question then is whether those instances where it holds
-- we have discussed Newton's theory of planetary motion, the
Schr\"odinger equation, biochemical reaction kinetics, the Hodgkin-Huxley equations or traffic dynamics -- are just lucky coincidences, or whether there are more systematic structural aspects of reality -- whatever that is -- that make our speculations sometimes so successful.

\section{The dream of mathematics}\label{dream}
Does this bring mathematics closer to its dream of a pure structural science, abstract and independent of specific content? We have argued already above that mathematics seeks regularities. Perhaps these regularities are discovered in specific domains, but hopefully,  they should apply also in other domains, and ideally, they should be universal. For instance, while Lie's  theory of symmetry groups was originally motivated by the symmetries of the systems of classical mechanics, the theory as such is abstract. In particular, through the work of W.Killing and E.Cartan, it lead to the classification of all possible continuous symmetries (see \cite{Haw} for the history). This turned out to be of fundamental importance for quantum mechanics, and even further for quantum field theory. It  permeates much of pure mathematics. It is also relevant in all applications where continuous symmetries arise. Even more generally, the notion of a group (Lie groups are  a special type of continuous groups), as developed by Gauss and Galois, first arose from studying solutions of algebraic equations, that is, a specific mathematical problem, but it then become something of a paradigm of a mathematical structure as such, and it is now important in almost all fields of mathematics, as well as in many areas of physics and other disciplines. 

Perhaps data science enables further steps in this direction. Here, even the original problems are no longer domain specific, as the symmetries and invariances of the systems of classical mechanics, but by their very nature already general and abstract when they apply to all kinds of large data sets. It is then natural that the mathematical tools developed to investigate these problems are at least as abstract as those problems themselves. Thus, data science may lead to an abstract theory of structures, that is, the purest form of mathematics. Incidentally, this brings fields of mathematics into focus that hitherto have been considered as applied and not pure mathematics. Statistics is an example. On the one hand, it can be mathematically treated in terms of the geometry of families of probability distributions, see e.g. \cite{ANH,AJLS}, and on the other hand, high-dimensional statistics and machine learning  might also become a science of distributions of geometric objects in high-dimensional spaces, again independent of any specific content. But there should be mathematical structures even more abstract and general than that.

\subsection*{Acknowledgements}
I am very grateful to my partners in the ZIF project, Philippe Blanchard, Martin Carrier, Andreas
Dress, Johannes Lenhard and Michael
R\"ockner, for their insightful comments and
helpful suggestions which have contributed to shaping and
sharpening the
ideas presented here. Martin Carrier and Johannes Lenhard also provided very useful comments on an earlier version of my text. Some of my own work discussed here was supported
by the ERC Advanced Grant FP7-267087 and the EU Strep ``MatheMACS''.

\end{document}